\documentclass{article}
\usepackage{amsfonts,amssymb,latexsym}
\makeatletter
\def\title#1{\gdef\@title{\uppercase{#1}}}

\def\maketitle{\par
\begingroup
   \def\thefootnote{\fnsymbol{footnote}}%
   \def\@makefnmark{\hbox
       to\z@{$\m@th^{\@thefnmark}$\hss}}%
   \newpage
   \global\@topnum\z@
   \@maketitle
   \thispagestyle{romp}\@thanks\markright{\RHD}\markboth{\LHD}{\RHD}%
\endgroup
 \setcounter{footnote}{0}%
 \let\maketitle\relax
 \let\@maketitle\relax
 \gdef\@thanks{}\gdef\@author{}\gdef\@title{}\let\thanks\relax}

\def\@maketitle{\newpage
 \null
 \vskip 4em
 \begin{center}
  {\normalsize\bf{\@title}\vskip20pt\par}
  {\normalsize\rm\begin{tabular}[t]{c}\@author\end{tabular}\vskip5pt\par}
 \end{center}\vskip2pt}

\def\abstract{\begin{list}{}{\leftmargin17mm\topsep0pt}%
\item\small\hspace*{3ex}}
\def\endabstract{\end{list}\vskip18pt\par}

\def\section{\@startsection {section}{1}{\z@}{3ex plus 1ex minus
.2ex}{1.2ex plus .2ex}{\normalsize\bf}}
\def\subsection{\@startsection{subsection}{2}{\z@}{3ex plus 1ex minus
.2ex}{1.2ex plus .2ex}{\normalsize\bf}}

\def\@listI{\leftmargin\leftmargini \parsep 4\p@ plus1\p@ minus\p@
\topsep 8\p@ plus1\p@ minus2\p@
\itemsep 0\p@ plus0\p@ minus0\p@}

\def\thebibliography#1{\section*{{\normalsize\rm REFERENCES}}
\small\rm\list
 {[\arabic{enumi}]}{\settowidth\labelwidth{[#1]}\leftmargin\labelwidth
 \advance\leftmargin\labelsep\usecounter{enumi}}
 \def\newblock{\hskip .11em plus .33em minus -.07em}
 \sloppy\clubpenalty4000\widowpenalty4000
 \sfcode`\.=1000\relax}

\@addtoreset{footnote}{section}

\def\ps@myheadings{\let\@mkboth\@gobbletwo
\def\@oddhead{\hbox{}\hfil{\footnotesize\rm\rightmark}\hfil
\normalsize\rm\thepage}\def\@oddfoot{}\def\@evenhead{\normalsize\rm
\thepage\hfil{\footnotesize\rm\leftmark}\hbox{}\hfil
}\def\@evenfoot{}\def\sectionmark##1{}\def\subsectionmark##1{}}

\font\rompF=cmmi6
\def\ps@romp{\let\@mkboth\@gobbletwo
  \def\@oddhead{\small\rm Vol.\ \Volume\ (\Year)\hfil
  {\rompF REP\kern.5ptO\kern-.6pt
RTS\hspace{1.2ex}O\kern-.3ptN\hspace{1.2ex}MATHEMATIC\kern-1pt
AL\hspace{1.2ex}PHY\kern.4ptSIC\kern-.5ptS}\hfil
     No.\ \Number}\def\@oddfoot{\rm\hfil[\thepage]
     \hfil}\def\@evenhead{}\let\@evenfoot\@oddfoot}

\newcounter{defthM}
\newenvironment{definition}[1]{\stepcounter{defthM}\trivlist
   \item[\hskip19pt{\sc #1~\arabic{defthM}.}]\rm\hskip3pt}{\endtrivlist}
\newenvironment{theorem}[1]{\stepcounter{defthM}\trivlist
   \item[\hskip19pt{\sc #1~\arabic{defthM}.}]\it\hskip3pt}{\endtrivlist}

\def\thefootnote{\arabic{footnote}}

\def\RHD{[Author and title]}
\def\LHD{[Author and title]}
\def\Year{199X}
\def\Volume{XX}
\def\Number{X}
\makeatother

\newcommand{\be}{\begin{equation}}
\newcommand{\ee}{\end{equation}}
\newcommand{\ba}{\hspace*{-5pt}\begin{array}}
\newcommand{\ea}{\end{array}}
\newcommand{\p}{\partial}
\newcommand{\ds}{\displaystyle}
\title{ On symmetries of KdV-like evolution equations}
\author{{\sc Artur Sergyeyev}\thanks{\hspace*{2mm}Current address (as of January 9, 2007):
Silesian University in Opava, Mathematical Institute,
Na~Rybn\'\i{}\v cku 1, 746\,01~Opava, Czech Republic.
E-mail: {\tt Artur.Sergyeyev@math.slu.cz}}\\
Institute of Mathematics of NAS of Ukraine,\\
Tereshchenkivs'ka Str. 3, 252004 Kyiv, Ukraine\\
e-mail: arthurser@imath.kiev.ua, arthur@apmat.freenet.kiev.ua}

\begin{document}
\maketitle
\begin{abstract}
The $x$-dependence of the symmetries of (1+1)-dimensional scalar
translationally invariant evolution equations is described. The
sufficient condition of (quasi)po\-ly\-nomiality in time $t$ of the
symmetries of evolution equations with constant separant is found.
The general form of time dependence of the symmetries of KdV-like
non-linearizable evolution equations is presented.
\end{abstract}

\section{Introduction}

It is well known that provided scalar (1+1)-dimensional evolution
equation (EE) with time-independent coefficients possesses the
infinite-dimensional commutative Lie
algebra of {\it time-independent} non-classical symmetries, it is
either linearizable or integrable via inverse scattering transform
\cite{s, o}. This algebra is usually constructed with usage of the recursion
operator \cite{o},
but
it may be also generated by the repeated commuting of
{\it mastersymmetry} with few time-independent symmetries \cite{fu}.
In its turn, to
possess the mastersymmetry, EE in question must have (at
least one) polynomial in time $t$ symmetry.
This fact is one of the main reasons of growing interest to the
study of whole algebra of {\it time-dependent} symmetries of
EEs \cite{mbcf,fu1,mawx}.
\looseness=-1

However, it is very difficult to describe this algebra
even in the simplest case of scalar (1+1)-dimensional EE.
To the best of author's knowledge, in
the class of scalar {\it nonlinear} EEs the complete
algebras of time-dependent local symmetries were found only for KdV
equation by Magadeev and Sokolov \cite{ms} and
for KdV and Burgers equations by Vinogradov et al. \cite{v}. In
\cite{v} there were also proved two no-go theorems, which show, when
the symmetries of third order KdV-like and second order Burgers-like
EEs are exhausted by Lie ones.
\looseness=-1
%

Surprisingly enough, using only the invariance of given EE
under $\p /\p t$ and $\p /\p x$ and some simple observations on
the explicit form of symmetries, we can show (vide Theorem 6
below) that
any symmetry of KdV-like non-linearizable EE as a function of $t$ is
a linear combination of quasipolynomials (i.e.\  of the products of
exponents by polynomials). Moreover, for the class of
EEs with constant separant our approach allows to find a simple sufficient
condition for such a situation to take place (Theorem 3), and
for the generic EE (\ref{eq: eveq}) it enables us to describe
the dependence of its symmetries on variable $x$ (Theorem 1).
\looseness=-1
\section{Some general properties of symmetries of evolution equations}
Consider the scalar $(1+1)$-dimensional EE
\begin{equation} \label{eq: eveq}
\partial u / \partial t = F(t, u, u_{1}, \dots, u_{n}), \quad n \geq 2,
\end{equation}
where $u_l=\partial ^{l} u /\partial x^{l}$,$l \in \mathbb{N}$,
$u_{0}\equiv u$, and its symmetries, i.e.\  the right hand sides $G$ of EEs
\begin{equation} \label{eq: sym}
 \partial u / \partial \tau =G(x,t,u, u_{1}, \dots, u _{k}),
\end{equation}
compatible with equation~(\ref{eq: eveq}).
The biggest number $k$ such that $\p G /\p u_{k} \neq 0$ is called
the order of symmetry and is denoted as $k={\rm ord }\; G$. If $G$ is
independent of $u, u_{1}, \dots$, then we assume that ${\rm ord }\; G=0$.
We set $S^{(k)}=\{ G\in S|{\rm ord}\; G \leq k \}$ and $S =
\bigcup_{j=0}^{\infty} S^{(k)}$.
\looseness=-2

For any sufficiently smooth function $h(x, t, u, u_1, \dots, u_r)$
introduce the quantities \cite{s}
$$
h_{*} = \sum\limits_{i=0}^{r} \p h /\p u_{i} D^{i} \: \mbox{and} \:
\nabla_{h} =\sum\limits_{j=0}^{\infty} D^{j}(h) \p /\p u_{j},\:
\mbox{where} \: D= {\ds \frac{\p}{\p x}} + \sum\limits_{i=0}^{\infty}
u_{i+1} \p/\p u_{i}.
$$
Now we can define the Lie bracket, which endows $S$ by the structure
of Lie algebra, as
$$
\{ h,r \} = h_{*} (r) - r_{*} (h) = \nabla_{r} (h) - \nabla_{h} (r).
$$
This definition
differs from the conventional one
\cite{s, o, ma} by the sign, but is more suitable for our purposes.
Note that $S^{(1)}$ is Lie subalgebra in $S$.

Equation (\ref{eq: sym}) is compatible with (\ref{eq: eveq}) if and only if
\be \label{eq: comp}
\p G /\p t = \{ F, G \}.
\ee
%
From (\ref{eq: comp}) one may easily derive \cite{o} that
\be \label{eq: cr3}
\p G_{*} /\p t \equiv (\p G/\p t)_{*} = \nabla_{G}
(F_{*}) - \nabla_{F}(G_{*}) + [F_{*}, G_{*}],
\ee
where $\nabla_{F}(G_{*}) \equiv \sum\limits_{i,j=0}^{\infty} D^{j}(F)
{\ds \frac{\p^{2} G}{\p u_{j} \p u_{i}}} D^{i}$ and similarly for
$\nabla_{G} (F_{*})$; $[\cdot,\cdot]$ stands for the usual commutator
of linear differential operators.

Let ${\rm ord}\; G=k$. Then equating the coefficients at $D^{l}$ on both sides
of (\ref{eq: cr3}) yields
\be \label{eq: cr4}
\ba{l}
{\ds \frac{\p^{2} G}{\p u_l \p t}} = \sum\limits_{m=0}^{n} D^{m} (G) {\ds
\frac{\p^{2} F}{\p u_{m} \p u_{l}}} - \sum\limits_{r=0}^{k} D^{r}
(F) {\ds \frac{\p^{2} G}{\p u_{r} \p u_{l}}} \\[4mm]
+ \sum\limits_{j=\max(0, l+1-n)}^{k} \sum\limits_{i=\max(l+1-j,0)}^{n}
\Bigl[ {\rm C}_{i}^{i+j-l} {\ds \frac{\p F}{\p u_{i}}} D^{i+j-l}
\Bigl( {\ds \frac{\p G}{\p u_{j}}} \Bigr) \\[4mm]
 - {\rm C}_{j}^{i+j-l} {\ds \frac{\p G}{\p u_{j}}} D^{i+j-l} \Bigl(
{\ds \frac{\p F}{\p u_{i}}}\Bigr) \Bigr],
 \quad l=0, \dots, n+k-1,
\ea
\ee
where $\ds {\rm C}_{q}^{p} = \frac{q!}{p! (q-p)!}$ and we assume that
$\ds \frac{1}{(-s)!} = 0$ for $s \in \mathbb{N}$.

If $k \geq 2$, one easily obtains from (\ref{eq: cr4}) with $l=n+k-1,
\dots, n+1$ the formulas
\be \label{eq: der1}
\p G /\p u_k = c_{k} (t) (\p F/\p u_{n})^{k/n},
\ee
\be \label{eq: derr}
\p G/\p u_{i} = c_{i} (t)(\p F/\p u_{n})^{i/n} +
\sum_{p=i+1}^{k}\sum_{q=0}^{\lbrack \frac{p-i}{n-1} \rbrack}
\chi_{pq} (t, x, u, \dots, u_k) {\ds \frac{\p^{q} c_{p}}{\p t^{q}}},
i=2, \dots,k-1,
\ee
where $c_{j} (t)$ are arbitrary functions of $t$ (cf. \cite{s,o}).

Furthermore, we see that by virtue of (\ref{eq: cr3})
\be \label{eq: mod}
[\nabla_{F} - F_{*}, G_{*}] =0\: \pmod{D^{p}}, \quad p= \max(k,n)
\ee
Equating the coefficients at powers of $D$ in (\ref{eq: mod}) yields the
equations (\ref{eq: cr4}) with $l=p+1, \dots, n+k-1$, from which we
may find $\p G/\p u_{i}$, $i=\max(k-n+1,2),\dots,k$. By (\ref{eq:
derr}), the only arbitrary elements, which they
may contain, are functions $c_{i} (t)$, while their dependence on
$x, u, u_1, \dots$ is {\it uniquely} determined from (\ref{eq: mod}).

On the other hand, since $F$ is $x$-independent, the
existence of the solution $F_{*}$ of (\ref{eq: mod}) with $p=n$ guarantees the
solvability of the equations for $\p G/\p u_{i}$, $i=\max(k-n+2,
2),\dots,k$, in terms of functions of $u, \dots, u_{k}$ and $t$ only
(cf.\ \cite{i}).
Therefore, $\p G/\p u_{i}$, $i=\max(k-n+2, 2),\dots,k$, are
$x$-independent.
In particular, any $G \in S^{(n)}$ has the form
\be \label{eq: low}
G = g(t,u, \dots, u_{k}) + \Phi (t,x,u,u_1).
\ee
\looseness=-1

Thus, for any symmetry $G \in S$
$\p G/\p x \in S$
and ${\rm ord}\; \p G/\p x \leq \max(1, {\rm ord}\; G - n + 1)$.
Applying this result to $\tilde G = \p G/\p x$ and so on, we
obtain that $\p^{r} G/\p x^{r} \in S^{(n)}$
and hence is of the form (\ref{eq: low}), if $r = r_{k,n,1}$, where
for $q=0,1$
$$
r_{k,n,q} = \left\{
\begin{array}{l}
\Bigl[{\ds \frac{k}{n-1}} \Bigr] \: \mbox{for} \:
k \not\equiv 0,\dots, q \!\pmod{n-1},\\
\max(0, \Bigl[ {\ds \frac{k}{n-1}}
 \Bigr] -1) \:
\mbox{for} \: k \equiv 0,\dots,q \!\pmod{n-1},
\end{array} \right.
\mbox{and}\: r_{k,n,-1} = \Bigl[{\ds \frac{k}{n-1}} \Bigr];
$$
$\lbrack s \rbrack$ denotes here the integer part
of the number $s$. The integration of
$\p^{r} G/\p x^{r}$ $r$ times with respect to $x$, taking into
account the above,
yields the following result:
\begin{theorem}{Theorem} \label{t1}
Any symmetry $G$ of order $k$ of (\ref{eq: eveq}) may be represented
in the form
\be \label{eq: rep}
G= \psi (t,x,u, u_1) + \sum_{j=0}^{s} x^{j} g_{j} (t,u, \dots, u_{k-
j(n-1)}), \quad s \leq r_{k,n,1}.
\ee
\end{theorem}
\begin{definition}{Remark} \label{rk1}
In complete analogy with the above, one may
show that if
\be \label{eq: deriv}
\p F/\p u_{n-i} = \phi_{i} (t),\quad i=0, \dots,j,
\ee
where $\phi_{i} (t)$ are arbitrary functions of $t$,
then in (\ref{eq: mod}) $p =\max(k, n-1-j)$ and it is possible to find
from (\ref{eq: mod}) $\p G/\p u_{i}$,
$i=\max(k-n+2, \max(1-j,0)),\dots,k$, which again turn out to be
$x$-independent, and hence (cf. \cite{i}) in (\ref{eq:
rep}) $s \leq r_{k,n, -\min(1,j)}$ and $\psi$ satisfies
\be \label{eq: der2}
\p \psi/\p u_{r} =0, \quad r=\max(1-j,0), \dots, 1.
\ee
Note that if (\ref{eq: deriv}) holds true, (\ref{eq: der1}),
(\ref{eq: derr}) hold for $k \geq \max(1-j,0)$, $i=\max(1-j,0),
\dots, k-1$.
\end{definition}
%
\section{Symmetries of the equations with constant separant}
Let us turn to the particular case, when $\p F/\p t=0$ and $F$
has the form
\be \label{eq: kdv}
F= u_{n} + f(u, \dots, u_{n-1}),
\ee
i.e.\  when $F$
has a constant separant, equal to unity \cite{i}. Note that
any $F$
with constant separant, different
from unity, may be reduced to the form (\ref{eq: kdv}) by rescaling
of time $t$. For the sake of brevity
we shall refer to EE
(\ref{eq: eveq}) with $F$
(\ref{eq: kdv}) as to EE
with constant separant. Let us also mention that if $\p F/\p t=0$,
then in (\ref{eq: derr}) $\p \chi_{pq}/\p t=0$.
\looseness=-1

Assume that ${\rm ord}\; G \equiv k > n-1$. Then, by (\ref{eq: kdv})
and (\ref{eq: der1}), (\ref{eq: cr4}) with $l=k$ reads
\be \label{eq: time1}
n D(\p G/\p u_{k-n+1}) = \p c_{k} (t)/\p t  + R,
\ee
where $R$ stands for the terms which depend only on $F$ and its
derivatives and on $\p G/\p u_{i}$, $i=k-n+2,\dots, k$.
Moreover, $R = D(K)$ for some $x$-independent $K$,
as it follows from the fact that if $F$ has a constant separant,
$F_{*}$ is solution of (\ref{eq: mod}) with $p=n-1$. Really, if the
term $\p G_{*}/\p t$ in (\ref{eq: cr3})
would be absent, $G_{*}$ would satisfy (\ref{eq: mod}) with $p=n-1$
and the equation for $\p G/\p u_{k-n+1}$ would
be solvable in terms of functions of $t, u, u_{1}, \dots$ (cf.\ the proof
of Theorem 1 and \cite{i}).
But the only term in (\ref{eq: time1}),
generated by $\p G_{*}/\p t$, is $ \p c_{k} (t)/\p t$, while $R$
is the same as if it would be in absence of $\p G_{*}/\p t$
in (\ref{eq: cr3}). Hence, $R =D(K)$, $\p K/\p x =0$, and
\be \label{eq: deriv1}
\p G/\p u_{k-n+1} = (x/n) \p c_{k} (t)/\p t + K  + c_{k-n+1} (t).
\ee

Since $\p G/\p u_{i}$,
$i=k-n+2,\dots,k$, are $x$-independent by Theorem 1, by (\ref{eq:
deriv1}) ${\rm ord}\; \p G/\p x = k-n+1$ and
$\p^{2} G/\p u_{k-n+1} \p x=  (1/n) \p c_{k} (t)/\p t$.
\looseness=-1

Iterating this process shows that for $r=r_{k,n,0}$
$Q=\p^{r} G/\p x^{r}\in S^{(n-1)}$ and
\be \label{eq: lead1}
\p Q/\p u_{q} = (1/n^{r}) \p^{r} c_{k} (t)/\p t^{r},\: q \equiv {\rm
ord}\; Q.
\ee
\looseness=-1
%
\begin{theorem}{Theorem} \label{t2}
If the symmetries from $S^{(n-1)}$ of the equation (\ref{eq:
eveq}) with constant separant either are all polynomial in $t$ or are
all linear combinations of quasipolynomials\footnote{We call
quasipolynomials the products $\exp (\lambda t) P(t)$, where $\lambda \in
\mathbb{C}$ and $P$ is a polynomial.} in $t$, then so does any
symmetry of this equation.
\end{theorem}
\looseness=-1
\noindent{\it Proof.} If the conditions of theorem are
fulfilled, then by (\ref{eq: lead1}) for any symmetry $G$, $k \equiv
{\rm ord}\; G \geq 1$, the function $c_{k} (t) \equiv \p G/\p u_{k}$
is either polynomial or linear combination of quasipolynomials in $t$.
Hence, there exists a differential operator $\Omega =
\sum_{l=0}^{m} a_{l} \p^{l} /\p t^{l}$,
$a_{l} \in \mathbb{C}$, such that $\Omega (c_{k}(t))=0$. In particular, if
$c_{k} (t)$ is polynomial in $t$ of order $p$, we may choose $\Omega_{0}
=\p^{p+1}/\p t^{p+1}$ as $\Omega$.
\looseness=-1

Now assume that the theorem is already proved for the symmetries from
$S^{(k-1)}$ (it is obviously true for $k \leq n$).
For EE (\ref{eq: eveq}) with $F$ (\ref{eq: kdv}) $S^{(k)}$ is closed under
$\p/\p t$,
and therefore $\Omega (G) \in S^{(k)}$.
Moreover, since $\Omega (c_{k} (t)) = 0$,
${\rm ord}\; \Omega (G) \leq k-1$ and hence, by our assumption, $\tilde
G \equiv \Omega (G)$ is either polynomial or linear combination of
quasipolynomials in $t$. Obviously, so does any solution $R$ of the equation
$\Omega (R) = \tilde G$, including $R=G$. If all
the elements of $S^{(n-1)}$ are polynomial in $t$, then so does
$c_{k} (t)$ and hence the polynomiality of $G$ in $t$ is guaranteed,
because we may take $\Omega = \Omega_0$ and because $\tilde G$ is
polynomial in $t$ by our assumption. The induction by $k$, starting
from $k=n$, completes the proof. $\square$
\looseness=-1

Theorem 3 is a natural generalization of the result of \cite{ms} on
polynomiality in $t$ of symmetries of KdV equation.
It gives a very simple sufficient condition for all the symmetries
of a given EE with constant separant to be polynomial
in time $t$. Note that in such a situation all the time-dependent
symmetries of EE in question may be constructed via the so-called generators
of degree $s$ for different $s \in \mathbb{N}$, using the results of
Fuchssteiner \cite{fu}.
\looseness=-1
\section{Symmetries of KdV-like equations}
Now let us consider
the equations with constant separant, whose $f$ satisfies
\be \label{eq: der3}
\p f/\p u_{n-1} ={\rm const}.
\ee
We shall call the EEs (\ref{eq: eveq}) with $F$ (\ref{eq:
kdv}), satisfying (\ref{eq: der3}),
{\it KdV-like}, since the famous Korteweg -- de Vries
equation has the form (\ref{eq: eveq}) with $F$ (\ref{eq: kdv}),
where $n=3$ and $f = 6 u u_{1}$ obviously satisfies (\ref{eq: der3}).
\looseness=-1

Let $G$ be the symmetry of KdV-like EE (\ref{eq: eveq}). Analyzing
the leading term of $\p^{r} G/\p x^{r}$, $r=\Bigl[ \frac{{\rm ord}\;
G}{n-1} \Bigr]$, like in the proof of Theorem 3, we obtain the
following statement:
\begin{theorem}{Corollary} \label{cor1}
If the symmetries from $S^{(n-2)}$ of KdV-like equation (\ref{eq:
eveq}) either are all polynomial or are all linear
combinations of quasipolynomials in $t$, then so does any symmetry of
this equation.
\looseness=-1
\end{theorem}
\begin{definition}{Example}
Consider third order formally integrable
nonlinear KdV-like EEs \cite{mss}:
\begin{displaymath}
\ba{l}
u_{t} = u_{3} + u u_{1},\\
u_{t} = u_{3} + u^{2}_{1} + c,\\
u_{t} = u_{3} + u^{2} u_{1} + c u_{1},\\
u_{t} = u_{3} + u^{3}_{1} + c u_{1} + d,\\
u_{t} = u_{3} - u^{3}_{1}/2 + (a \exp( 2 u) + b \exp (- 2 u) + d) u_{1},
\ea
\end{displaymath}
where $a,b,c,d \in \mathbb{C}$. All the symmetries of these EEs
are polynomial in $t$ by Corollary 4, since so do their
symmetries of orders 0 and 1.
\end{definition}

Now let us analyze in more detail the general form of time dependence
of symmetries of KdV-like EE (\ref{eq: eveq}).
Assume that the EE in question may not be linearized by means of
contact transformations (for the sake of brevity we shall call it
{\it non-linearizable}). Then ${\rm dim}\; \Phi \leq n$ \cite{ma},
where $\Phi = \{ \varphi(x,t) | \varphi(x,t) \in S \}$. Let us show
that in such a case ${\rm dim}\; S^{(k)} < \infty$ for any $k=0,1,2,\dots$.
\looseness=-1

Let $G \in S^{(k)}/S^{(k-1)}$, $k > n-2$. Then, obviously, it
is completely determined by its leading term $h_{k} (t)\equiv \p G/\p
u_{k}$. Like the above, but taking into account
Remark 2, we may show that for KdV-like EE (\ref{eq: eveq})
$Q = \p^{r} G/\p x^{r} \in S^{(n-2)}$,
if $r=\Bigl[ \frac{k}{n-1} \Bigr]$, and
\be \label{eq: lead2}
\p Q/\p u_{q} \equiv c_{q} (t) = (1/n^{r}) \p^{r} h_{k} (t)/\p
t^{r}, q \equiv {\rm ord}\; Q.
\ee
Since $h_{k}$ satisfies (\ref{eq: lead2}), for $k > n-2$
\be \label{eq: dim0}
{\rm dim}\; S^{(k)}/S^{(k-1)} \leq {\rm dim}\; S^{(k_0)} + \Bigl[
\frac{k}{n-1} \Bigr],
k_0=k- \Bigl[ {\ds \frac{k}{n-1}} \Bigr] (n-1),
\ee
whence ${\rm dim}\; S^{(k)} < \infty$ for $k=0, \dots, n-2$
implies the same result for any $k$.
\looseness=-1

By Theorem 1 and Remark 2 for KdV-like EE (\ref{eq:
eveq}) (\ref{eq: der1}) and (\ref{eq: derr}) for $k \leq n-2$ read
\begin{eqnarray}
\label{eq: simple}
\p G/\p u_{i} = c_{i} (t) + \sum_{p=i+1}^{k}
\chi_{p} (u, \dots, u_k) c_{p} (t), \quad i=0, \dots,k-1, \\
\p G/\p u_{k} = c_{k} (t),
\end{eqnarray}
and thus any symmetry $G$ of order $k \leq n-2$ has the form
\be \label{eq: rep2}
G= \psi (t,x) + g_{0} (t,u, \dots, u_{k}).
\ee
Without loss of generality we can assume that the function $g_{0}$ is
completely determined by $\p G/\p u_{i}$, $i=0, \dots, k$. Since $\p
\psi/\p x \in \Phi$, we have
\be \label{eq: psi}
\psi(t,x) = \gamma(t) + \sum_{p=1}^{{\rm dim }\; \Phi} a_{p} \int_{0}^{x} dy\;
\varphi_{p}(y,t),
\ee
where $a_{p} \in \mathbb{C}$, $\gamma
(t)$ is arbitrary function of $t$, and $\varphi_{p}(x,t)$, $p=1, \dots, {\rm
dim}\; \Phi$, stand for some basis in $\Phi$.
\looseness=-1

The substitution of $G$ (\ref{eq: rep2}) with $\psi$ (\ref{eq:
psi}) into equations (\ref{eq: cr4}) with $l=0, \dots, n-1$
and into (\ref{eq: comp}) yields in final account the system of first
order linear ordinary differential equations in $t$ (and, possibly,
algebraic equations) for $c_{i} (t)$, $i=0, \dots, k$ and
$\gamma (t)$. Note that we must use (\ref{eq: comp}) in order to
obtain an ODE of the form $\p \gamma (t)/\p t =
\cdots$, allowing to find $\gamma(t)$.

Obviously, the general solution of this system of ODEs for $k \leq
n-2$ may contain at most $N_{k,n}= {\rm dim}\; \Phi + k + 2$
arbitrary constants (including $a_{p}$, $p=1, \dots, {\rm dim}\; \Phi$).

Hence, ${\rm dim}\;S^{(k)} \leq N_{k,n} < \infty$ for $k \leq n-2$ and
thus by (\ref{eq: dim0}) for any $k=0,1, \dots$
\be \label{eq; dim2}
 {\rm dim}\; S^{(k)} =
{\rm dim}\; S^{(k_0)} + \sum_{j=k_0+1}^{k} {\rm dim}\;
S^{(j)}/S^{(j-1)} < \infty.
\ee

Thus, for any $k$ the space $S^{(k)}$ is
finite-dimensional. Since in addition this space is
invariant under $\p /\p t$, the dependence of its elements on $t$
is completely described by Theorem 3.1 \cite{ss}. Namely,
any symmetry of order $k$ of KdV-like non-linearizable EE (\ref{eq:
eveq}) is a linear combination of ${\rm dim}\; S^{(k)}$ linearly
independent symmetries of the form
\be \label{eq: time}
H = \exp(\lambda t) \sum_{j=0}^{m} t^{j} h_{j} (x, u, \dots,
u_{k}),
\lambda \in {\mathbb{C}}, m \leq {\rm dim}\; S^{(k)}-1.
\ee
%
\begin{theorem}{Theorem} \label{t3}
 For any non-linearizable KdV-like EE (\ref{eq: eveq})
${\rm dim}\; S^{(k)}
< \infty , k=0,1,2,\dots$
and any symmetry $Q$ of order $k$
is a linear
combination of the symmetries (\ref{eq: time}).
\looseness=-1
\end{theorem}
Thus, all the symmetries of non-linearizable KdV-like EEs are linear
combinations of quasipolynomials in $t$. This partially recovers the result of
Corollary 4, but this corollary still remains of interest, providing
the convenient sufficient condition of {\it polynomiality} of symmetries in
time $t$.
\looseness=-1

It is interesting to note that some general
properties of time-dependent symmetries, which are linear
combinations of the expressions (\ref{eq: time}),
were studied by Ma \cite{mawx}. However, while he considered this
form as given {\it a priori}, we have {\it proved} that all the
symmetries of non-linearizable KdV-like EE (\ref{eq: eveq})
indeed have this form.

Moreover, acting on any symmetry (\ref{eq: time}) by $(\p /\p
t - \lambda)^{m}$ for $\lambda \neq 0$ or by $\p^{m-1} /\p t^{m-1}$
for $\lambda = 0$, we obtain the symmetry which is either linear or
exponential in $t$. Hence, there is a very simple test of
existence of {\it any}
time-dependent symmetries for given non-linearizable KdV-like EE
(\ref{eq: eveq}). Namely, it suffices to check whether there exist the
symmetries of the form
\be \label{eq: case1}
G= \exp (\lambda t) Q_0 , \lambda \in {\mathbb{C}}, \lambda \neq 0
\ee
or of the form
\be
\label{eq: case2} G = G_0 + t G_1, G_1 \neq 0,
\ee
where $Q_0$, $G_0$ and $G_1$ are time-independent. If
non-linearizable KdV-like EE
(\ref{eq: eveq}) (with time-independent coefficients!) has no
time-dependent symmetries of the form (\ref{eq: case1}) or (\ref{eq:
case2}), then it has no time-dependent symmetries at all (but of course it
may have time-independent symmetries).
\looseness=-1

The substitution of (\ref{eq: case1}) and (\ref{eq: case2}) into
(\ref{eq: comp}) yields
\begin{eqnarray} \label{eq: scale}
\{ F, Q_0 \} =\lambda Q_0,\\
\label{eq: mast}
\{F, G_0 \} = G_1, \quad \{ F, G_1 \} =0.
\end{eqnarray}

In the first case $F$ is called scaling symmetry (or conformal
invariance \cite{oew}) of $Q_0$. However, known
scaling symmetries $F$ of integrable hierarchies, such as
KdV, depend usually only on $x,u,u_1$ but not on $u_2$ and higher
derivatives \cite{oew} and hence do not generate EEs of the form
(\ref{eq: eveq}), which we consider
here. We guess that if KdV-like EE (\ref{eq: eveq}) is
non-linearizable and integrable, there exist no functions $Q_0$, which
satisfy (\ref{eq: scale}) with $\lambda \neq 0$. Moreover, it is believed
\cite{mbcf} that in such a case the only polynomial in $t$ symmetries
(\ref{eq: sym}) that EE (\ref{eq: eveq}) may possess are those
linear in $t$.
\looseness=-1

Now let us consider the second case. Assume that there exists some
commutative
algebra $Alg$ of time-independent symmetries of KdV-like
non-linearizable EE (\ref{eq: eveq}), such that for any $K \in Alg$
the Lie bracket $\{G_0,K \} \in Alg$.
Then $G_0$ is mastersymmetry of (\ref{eq: eveq}), and hence (\ref{eq:
eveq}) possesses (under some extra conditions, vide \cite{fu}) the infinite
set of time-independent symmetries and is probable to be integrable
via inverse scattering transform. Let us mention that the condition of
commutativity of $Alg$ may be rejected if $G_0$ is scaling symmetry
of $F$, i.e.\  $\{F, G_0 \} = \mu F$ for some $\mu \in {\mathbb{C}}$,
$\mu \neq 0$ \cite{oew}.
\looseness=-1

\begin{theorem}{Conjecture}
For any KdV-like non-linearizable evolution
equation (\ref{eq: eveq}) either all its symmetries are polynomial in
$t$ or all they are linear combinations of exponents in $t$.
\looseness=-1
\end{theorem}
\section*{Acknowledgements}
It is my pleasure to express deep gratitude to
Profs.~A.G.~Nikitin and R.Z.~Zhdanov and Dr.~R.G.~Smirnov for
the fruitful discussions on the subject of this work.
I would also like to thank the organizers of XXX Symposium on
Mathematical Physics, where this work was presented, for their
hospitality.

\end{document}